\documentclass[10pt]{article}
\usepackage[margin=1in]{geometry}
\usepackage{amsmath,amssymb,amsfonts}
\usepackage{graphicx}
\usepackage{booktabs}
\usepackage{siunitx}
\usepackage{algorithm}
\usepackage{algpseudocode}
\usepackage{cite}
\usepackage{hyperref}
\usepackage{xcolor}
\usepackage{csvsimple}
\usepackage{caption}
\usepackage{subcaption}
\usepackage{enumitem}
\usepackage{authblk}
\usepackage{float}
\usepackage{mathrsfs}
\usepackage{longtable}
\usepackage{amsmath,amssymb,amsthm}

\setlist{nosep}

\title{Adaptive $\mathcal{C}h$ Method with Local Coupled Multiquadrics for Solving Partial Differential Equations}
\author{ Ahmed E. Seleit}
\date{ }
\affil{Mechanical Engineering, Johns Hopkins University}
\graphicspath{{aml-covers-vs-shape/figs/}}

\theoremstyle{remark}

\begin{document}
\maketitle

\begin{abstract}
We present a new adaptive collocation scheme for solving partial differential equations based on Local Coupled Multiquadrics (LCMQs) within a covers-and-nodes framework. The method, referred to as the Adaptive $\mathcal{C}h$ Method, automatically prioritizes adjusting the local cover size $\mathcal{C}$ then refines local nodal spacing $h$ to achieve a prescribed tolerance. Numerical examples for one- and two-dimensional Poisson problems demonstrate accurate solutions across a wide range of shape parameter values, while preserving the advantages of local collocation. The proposed approximation approach is truly meshless, requiring no element, connectivity or continuity to construct trial functions or weights.
\end{abstract}

\section{Introduction}

Radial basis functions (RBFs) are smooth, radially symmetric functions originally developed for multivariate interpolation \cite{buhmann2000radial}. Their meshfree nature makes them especially attractive for high-dimensional problems and irregular geometries, since approximation can be constructed without elements, grids, or connectivity. Building on these properties, Kansa introduced the unsymmetric collocation method, which applied RBFs directly to the discretization of partial differential equations (PDEs) \cite{kansa1990multiquadrics, kansa1990scattered}. This approach opened a versatile framework for solving boundary-value and evolution problems using scattered nodes alone. Since then, numerous extensions have been developed, including local variants such as RBF-FD and partition-of-unity formulations to alleviate the ill-conditioning of global systems while preserving high-order accuracy; stability for hyperbolic problems has been analyzed with discrete $\ell_2$ estimates and oversampling strategies~\cite{Tominec2025}. Additional advances include localization and multi-symplectic schemes for long-time conservative dynamics~\cite{Nikan2021,Zhang2022}, meshless space–time backward substitution for advection–diffusion~\cite{Lin2021}, and fast solvers for nonlinear radiation diffusion~\cite{Wang2023}. Applications span elastoplastic analysis of functionally graded materials and thermo-mechanical stability of composites~\cite{Liu2022,Shukla2022}, capillary surface computation via pseudo-spectral methods~\cite{Treinen2023}, meshless cubature and moment computation on polygons~\cite{Sommariva2021}, and reduced-quaternionic function frameworks relevant to time-dependent wave equations~\cite{Morais2023}. Lee \emph{et al.} \cite{lee2003local} developed a \emph{local} multiquadric method in which a global banded matrix assembled on overlapping covers, thereby avoiding large dense RBF gramians. Collectively, these developments establish RBF collocation and its localized variants as flexible, accurate, and broadly applicable solvers for PDEs.\\

A central practical issue in RBF approximation is selecting the shape parameter, which mediates a trade-off between approximation accuracy (small parameter, flatter basis) and numerical conditioning (large parameter)~\cite{Esmaeilbeigi2014}. Variable and locally adaptive strategies have been proposed to tune this parameter within Kansa-type discretizations for nonlinear boundary-value problems, achieving significant time savings via analytic Jacobians in the nonlinear solve~\cite{Karageorghis2022}. Optimization-driven selection for local RBF collocation employs surrogate error objectives and particle-swarm search to assign spatially varying parameters that track solution features~\cite{Li2023}. Interval-based selection rules avoid delicate error-norm minimization by exploiting practical convergence behavior and stability considerations~\cite{Biazar2017}. Spectral and conditioning viewpoints connect the parameter to interpolation matrix properties and error behavior, yielding criteria for choosing both the kernel type and its parameter~\cite{Haq2018,Chen2023}. Alternative metaheuristics, such as genetic algorithms combined with regularized linear solves, offer additional pathways to robust choices~\cite{Esmaeilbeigi2014}. Recent progress on variable-projection formulations for general RBF neural networks also informs parameter learning in separable nonlinear least-squares settings~\cite{Zheng2023}. Another method was developed to find the a robust shape parameter for RBF approximation by: first, constructing an optimization problem to find a value of the shape parameter that leads to a bounded condition number then employ a data driven method that controls the condition number \cite{Maria2025}. Overall, these contributions provide theory-guided and data-driven tools that make shape selection systematic rather than ad hoc.\\

Zhang \cite{zhang2019accurate} introduced \emph{Coupled} RBFs (CRBFs) that reduce sensitivity to the shape parameter achieving accuracy and stable conditioning. Moreover, the author highlighted that Coupled Multiquadrics (CMQs) demonstrated better stability and performance compared to its counterparts. In contrast to standard RBFs that demand careful tuning of the shape parameter to avoid instability or poor convergence\cite{fornberg2004stable,fornberg2013stable,wright2017stable,cheng2012multiquadric,chen2018sample,rippa1999algorithm,ling2009subspace,sarra2009multiquadric,wang2018hausdorff,liu2018optimal,cheng2000drbem}, CRBFs remove this burden \cite{seleit2024shape, seleit2022shape}.\\

Foundational work introduced hybrid global/compact schemes with adaptive node refinement and support selection to balance stability and cost for evolving PDE solutions~\cite{Ling2012}. Adaptive meshfree integration via mapping techniques enabled accurate weak-form evaluations on complex domains without background cells, making several formulations truly meshless~\cite{Racz2012}. H-adaptive local collocation with simple but effective error indicators demonstrated substantial accuracy–performance gains for transport and Burgers-type problems~\cite{Kosec2011}. Subsequent procedures blended superposition-based error comparison and residual subsampling to add or remove collocation points iteratively~\cite{Cavoretto2020}. Early adaptive radial-basis algorithms for time-dependent PDEs further exploited the grid-free nature to place points where steep features develop~\cite{Sarra2005}. Beyond forward problems, adaptive mesh-free discretizations have been shown effective for stabilizing ill-conditioned inverse problems such as gravity inversion through modified RBFs and hybrid bases~\cite{Liu2023}. In sum, adaptive RBF frameworks consistently reduce computational complexity while preserving (and often improving) accuracy on complex, multiscale, or moving-boundary problems.\\

In this work, we introduce a novel approach that integrates: (i) the stability and low shape parameter sensitivity of CMQs \cite{zhang2019accurate}, with (ii) the idea of using local covers to solve PDEs \cite{lee2003local}, and (iii) adaptivity inspired by the \(ph\)-refinement method developed for optimal control \cite{patterson2015ph}. We automatically prioritize the adjustment of the cover size $\mathcal{C}$ then refine the nodal distance $h$ to accurately solve PDEs. We refer to this integrated approach as the Adaptive $\mathcal{C}h$ method. Our formulation avoids Partition of Unity weight terms in the differential operator \cite{cavoretto2019adaptive}, and assembles a banded global system with bandwidth controlled by the maximum cover size, yielding a simple implementation and predictable sparsity.

\section{LCMQ Covers}
Let $\Omega\subset\mathbb{R}^d$ be a bounded domain with boundary $\partial\Omega$ and $X=\{x_i\}_{i=1}^N\subset\Omega$ be distinct nodes. We seek $u:{\Omega}\to\mathbb{R}$ that satisfies
\begin{equation}
\mathcal{L}u(x)=f(x),\quad x\in\Omega,\qquad \mathcal{B}u(x)=g(x),\quad x\in\partial\Omega,
\end{equation}
where $\mathcal{L}$ is a linear differential operator and $\mathcal{B}$ collects Dirichlet/Neumann boundary operators. We denote by $u_h$ the meshless local–CRBF approximant introduced below.

For each node $x_i$, we associate a finite \emph{cover} $\Omega_{x_i}= \{x_{i_j}\}_{j=1}^{\mathcal{C}_i}\subset X$ of cardinality $\mathcal{C}_i$, typically chosen as the $\mathcal{C}_i$ nearest neighbors of $x_i$. Let $\rho_i:=\max_{x_{i_j}\in\Omega_{x_i}}\|x_{i_j}-x_i\|$ and $B_i:=\{x\in\Omega:\|x-x_i\|\le \rho_i\}$. The collection of covers is \emph{admissible} if (i) the union of cover regions embraces the domain and (ii) each target point involved in collocation belongs to at least one cover:
\begin{equation}
\Omega \;=\; \bigcup_{i=1}^N B_i,\qquad
\,x\in\Omega\;\iff\exists\,i:\;x\in B_i.
\end{equation}

\subsection{Local support}
The local approximation at a target point $x$ 
uses only the nodal values from $\Omega_{x_i}$, i.e.\ its support is restricted 
to the chosen cover. Note that covers generally overlap: a node $x_j$ may belong to 
several covers $\Omega_{x_i}$, leading to a sparse but 
globally coupled system. Let $\phi(r;c)$ denote a CRBF with shape parameter $c>0$, 
where $r = \|x - x_{i_j}\|$ is the distance between $x$ and a cover node 
$x_{i_j} \in \Omega_{x_i}$. In this work, without loss of generality \cite{zhang2024novel}, we adopt the coupled multiquadric (CMQ)~\cite{zhang2019accurate},
\begin{equation}
\phi_{\mathrm{CMQ}}(r;c) \;=\; \sqrt{1+\big(r/c\big)^2}\;+\; r^5.
\end{equation}
Let $U_i\in\mathbb{R}^{\mathcal{C}_i\times \mathcal{C}_i}$ be the local kernel matrix, $(U_i)_{j\ell}=\phi_{\mathrm{CMQ}}(\|x_{i_j}-x_{i_\ell}\|;c)$, where $j,\ell=1,\dots,\mathcal{C}_i$. Let
\begin{equation}
\psi_i(x)\;=\;\big[\phi_{\mathrm{CMQ}}(\|x-x_{i_1}\|;c),\ldots,\phi_{\mathrm{CMQ}}(\|x-x_{i_{\mathcal{C}_i}}\|;c)\big]^\top.
\end{equation}
The local shape vector is $W_i(x)=U_i^{-1}\psi_i(x)$ and the approximation (restricted to the union of covers) is
\begin{equation}
u_h(x) \;=\; \sum_{j=1}^{\mathcal{C}_i} W_{i,j}(x)\,u(x_{i_j}).
\end{equation}
At interior collocation points $x_i$, we enforce $\mathcal{L}u=f$ by
\begin{equation}
\sum_{j\in\Omega_{x_i}} \big(\mathcal{L}W_{i,j}\big)(x_i)\,u(x_j) \;=\; f(x_i).
\end{equation}

Let $\mathcal{L}$ be any second-order linear differential operator acting on $x$ (e.g.\ the Laplacian). Since $U_i$ is independent of $x$, we have $(\mathcal{L}W_i)(x)=U_i^{-1}(\mathcal{L}\psi_i)(x)$. Hence,
\begin{equation}
\sum_{j\in\Omega_{x_i}} \big(\mathcal{L}W_{i,j}\big)(x_i)\,u(x_j)
= \sum_{j=1}^{\mathcal{C}_i} \Big[\,U_i^{-1}\,(\mathcal{L}\psi_i)(x_i)\,\Big]_j\,u(x_{i,j}) \;=\; f(x_i),
\end{equation}
\noindent where $(\mathcal{L}\psi_i)(x)$ is defined component-wise by $\big((\mathcal{L}\psi_i)(x)\big)_j=\mathcal{L}\!\big[\phi(\|x-x_{i,j}\|;c)\big]$.
Boundary data are imposed directly. This procedure yields a sparse, banded global system.

\section{Adaptive \texorpdfstring{$\mathcal{C}h$}{} Method on LCMQ Covers}
We present the Adaptive $\mathcal{C}h$ method on a generic domain $\Omega\subset\mathbb{R}^d$. We employ a simple logic: when the a~posteriori indicator is unacceptable, first attempt $\mathcal{C}$-enrichment by increasing local cover cardinality; if insufficient, perform $h$-refinement by inserting nodes.

Let $\{K_g\}_{g=1}^{G}$ be a nonoverlapping partition of $\Omega$ into \emph{marker cells} used only for error estimation and node insertion. Denote the cell diameter by $h_g:=\operatorname{diam}(K_g)$ and its center by $\zeta_g$. For a baseline cover size $\mathcal{C}$, define an \emph{enriched} solution $u_h^{\mathrm{enr}}$ obtained by using cover size $\mathcal{C}+\Delta p$, where $\Delta p$ is a small integer, when evaluating the residual. The error indicator on $K_g$ is
\begin{equation}
\mathscr{E}_g \;=\; h_g^2\,\big|(\mathcal{L}u_h^{\mathrm{enr}})(\zeta_g)-f(\zeta_g)\big|.
\end{equation}
We also permit a classifier
\begin{equation}
s_g \;=\; h_g^2\,\big|(\mathcal{L}u_h^{\mathrm{enr}}-\mathcal{L}u_h^{\mathrm{base}})(\zeta_g)\big|,
\end{equation}
to bias the decision toward $\mathcal{C}$-enrichement or $h$-refinement.
This is a \emph{truly meshless method} because our approximation spaces are constructed from node clouds via local CMQ kernels; no element mesh or connectivity/continuity is used to define trial functions or weights. $K_g$ is employed only for error marking and node insertion,
which is an auxiliary device and does not affect the meshless construction. The method and adaptive cycle are further explained in the following subsections and Algorithm \ref{alg:ACMQ}.

\begin{algorithm}[H] 
\caption{Generic-$d$ Adaptive \texorpdfstring{$\mathcal{C}h$}{Ch} LCMQ  Collocation}
\label{alg:ACMQ}
\begin{algorithmic}[1]
\Require domain $\Omega\subset\mathbb{R}^d$; initial nodes $X^{(0)}$; CMQ shape $c$; cover cardinality bounds $\mathcal{C}{\min},\mathcal{C}{\max}$; tolerance $\tau$; max splits per cycle $M$; fraction of bad cells trigger $\theta\in(0,1)$; enrichment step $\Delta p$; classifier threshold $\rho\in(0,1)$; small $\varepsilon>0$; adaptive cycle $k=0,1,2,\dots,\mathscr{K}$.
\State $k\gets 0$; \textbf{initialize per-cover sizes} $\mathcal{C}_i\gets \mathcal{C}{\min}$ for all nodes $x_i\in X^{(0)}$.
\Repeat
  \State Assemble and solve the local--CRBF system on $X^{(k)}$ with \textbf{per-node cover sizes} $\{\mathcal{C}_i\}$ to obtain $u_h^{(k)}$ \cite{lee2003local}.
  \State For each marker cell $K_g$, compute $\mathscr{E}_g$ and $s_g$ at its center $\zeta_g$ using \textbf{enriched per-node sizes} $\{\mathcal{C}_i+\Delta p\}$ \emph{for evaluation only}.
  \State Let $\mathcal{G}\gets\{g:\mathscr{E}_g>\tau\}$ (marked cells).
  \If{$\dfrac{|\mathcal{G}|}{G}\ge\theta$ \textbf{and} $(\exists\, i:\, \mathcal{C}_i<\mathcal{C}{\max})$}
    \Comment{local $\mathcal{C}$-enrichement}
      \ForAll{$K_g\in\mathcal{G}$}
        \If{$s_g \;\ge\; \rho\big(\mathscr{E}_g+\varepsilon\big)$}
          \ForAll{nodes $x_i\in K_g$}
            \State $\mathcal{C}_i \gets \min\{\,\mathcal{C}_i+\Delta p,\; \mathcal{C}{\max}\,\}$
          \EndFor
        \EndIf
      \EndFor
      \If{$\exists K_g\in\mathcal{G}\text{ with } s_g < \rho(\mathscr{E}_g+\varepsilon)$}
        \State Select up to $M$ cells in $\mathcal{G}$ with $s_g < \rho(\mathscr{E}_g+\varepsilon)$ and largest $\mathscr{E}_g$; insert centers:
        \State \hspace{1.2em}$X^{(k+1)}\gets X^{(k)}\cup\{\text{selected }\zeta_g\}$.
        \State Initialize $\mathcal{C}_i\gets \mathcal{C}{\min}$ for each newly inserted node; keep existing $\{\mathcal{C}_i\}$ unchanged.
      \Else
        \State $X^{(k+1)}\gets X^{(k)}$
      \EndIf
  \Else
      \Comment{local $h$-refinement}
      \State Select up to $M$ cells in $\mathcal{G}$ with largest $\mathscr{E}_g$ and insert their centers:
      \State \hspace{1.2em}$X^{(k+1)}\gets X^{(k)}\cup\{\text{selected }\zeta_g\}$.
      \State Initialize $\mathcal{C}_i\gets \mathcal{C}{\min}$ for each newly inserted node; keep existing $\{\mathcal{C}_i\}$ unchanged.
  \EndIf
  \State Enforce caps on $|X^{(k+1)}|$ and on all $\mathcal{C}_i\in[\mathcal{C}{\min},\mathcal{C}{\max}]$.
  \State $k\gets k+1$.
\Until{$\max_g \mathscr{E}_g\le\tau$ or $k=\mathscr{K}$.}
\end{algorithmic}
\end{algorithm}
\subsection{Algorithmic rationale and flow.}
Algorithm~\ref{alg:ACMQ} advances the solution on a fixed domain $\Omega$ by alternating \emph{local cover enrichment} ($\mathcal{C}$-refinement) and \emph{node insertion} ($h$-refinement) until a target a~posteriori tolerance is met. The core idea is to modify the local approximation power only where the PDE operator appears under-resolved, and to increase sampling density only where enrichment alone is unlikely to cure the error. To that end, each adaptive cycle consists of four conceptually simple stages: solve, estimate, decide, and update.

\subsubsection{Initialization (Alg.\ 1, lines 1--2).}
We start from an initial node set $X^{(0)}$ and assign to each node $x_i$ a local cover cardinality $\mathcal{C}_i$ (size of the local stencil) set to the minimum admissible value $\mathcal{C}_{\min}$. Using per-node $\mathcal{C}_i$ rather than a single global size is essential: it allows the method to spend degrees of freedom only where needed while preserving sparsity and stability of the local LCMQ collocation \cite{lee2003local}. The bounds $[\mathcal{C}_{\min},\mathcal{C}_{\max}]$ and a cap on the total number of nodes control computational cost.

\subsubsection{Solve with current covers (line 4).}
Given the current node cloud $X^{(k)}$ and per-node cover sizes $\{\mathcal{C}_i\}$, we assemble the local CRBF (LCMQ) collocation system and compute a discrete solution $u_h^{(k)}$. Because each row uses only its local cover, the global matrix remains sparse/banded and well-conditioned relative to global RBF collocation, while the Kronecker-delta property of the local shape functions simplifies the enforcement of essential boundary conditions.

\subsubsection{A classifier to choose between $\mathcal{C}$-Enrichment and $h$-refinement (line 7).}
As mentioned above, alongside $\mathscr{E}_g$ we compute
\[
s_g \;=\; h_g^2\,\left|(\mathcal{L}u_h^{\mathrm{enr}}-\mathcal{L}u_h^{\mathrm{base}})(\zeta_g)\right|.
\]
$s_g$ measures the sensitivity of the operator evaluation to local cover enrichment alone. If $s_g$ is comparable to (or a significant fraction of) $\mathscr{E}_g$, then much of the observed discrepancy may stem from insufficient local approximation power, suggesting $\mathcal{C}$-enrichment is the right remedy. Conversely, if $s_g$ is small while $\mathscr{E}_g$ is large, local enrichment is unlikely to help and new samples (nodes) are warranted. The small $\varepsilon$ in the test $s_g \ge \rho(\mathscr{E}_g+\varepsilon)$ prevents division-by-zero logic and stabilizes decisions when both numbers are extremely small.

\subsubsection{Choosing the dominant action at cycle level (lines 8--9 and 18).}
We next examine the \emph{fraction} of failing cells, $|\mathcal{G}|/G$ with $\mathcal{G}=\{g:\mathscr{E}_g>\tau\}$. If a significant portion of the domain violates the tolerance (i.e., $|\mathcal{G}|/G\ge \theta$) and there is headroom to enlarge some covers ($\exists\,i:\mathcal{C}_i<\mathcal{C}_{\max}$), we prioritize $\mathcal{C}$-enrichment. A single solve after moderate local cover increases can correct broad, smooth under-resolution without proliferating nodes. If, instead, only a small portion fails or all covers are already at the cap, we default to $h$-refinement.

\subsubsection{Local $\mathcal{C}$-enrichment and selective $h$-refinement (lines 10--17).}
Within the $\mathcal{C}$-branch, we enrich \emph{only} those marked cells whose classifier supports it: for each $K_g\in\mathcal{G}$ satisfying $s_g\ge \rho(\mathscr{E}_g+\varepsilon)$, we increase the cover sizes $\mathcal{C}_i$ for all nodes $x_i$ lying in $K_g$ by $\Delta p$ (clamped to $\mathcal{C}_{\max}$). This keeps the action local and proportional to need. However, if some marked cells have $s_g< \rho(\mathscr{E}_g+\varepsilon)$, we also perform targeted $h$-refinement: among those, we pick up to $M$ worst cells (largest $\mathscr{E}_g$) and insert their centers $\zeta_g$ as new nodes. The cap $M$ avoids excessive growth in a single cycle and reduces the risk of clustered nodes that degrade conditioning. New nodes start with the safest choice $\mathcal{C}_{\min}$, while existing $\mathcal{C}_i$ remain unchanged—this conservative initialization helps maintain a stable linear system and lets the indicator decide future enrichment for the newcomers.

\subsubsection{$h$-refinement branch (lines 19--23).}
If global enrichment is not indicated (few bad cells, or covers saturated), we insert up to $M$ centers of the worst marked cells and keep all existing covers unchanged. This mirrors classical adaptive $h$-refinement: sampling density is increased where the residual is largest, which directly reduces local discretization error and improves the conditioning of nearby local solves.

\subsubsection{Enforcing caps, limits, and stopping (lines 24--27).}
Each cycle enforces global budgets: a cap on the number of nodes and on all cover sizes keeps complexity under control. The algorithm terminates when the maximum indicator falls below $\tau$ or when the allotted number of cycles $\mathscr{K}$ is reached. In practice, these safeguards balance accuracy with computational cost: cover enlargement is cheap (no new unknowns) but eventually saturates. Alternating the two, with the simple classifier and fraction trigger, yields a robust and efficient path to the target tolerance.

\section{Numerical Experiments}
We assess the accuracy and robustness of the proposed method on canonical Poisson problems in one and two dimensions. Throughout, $u$ denotes the exact solution and $u_h$ the numerical approximation obtained on the current node set and covers. To avoid nodal bias, errors are evaluated on an independent uniform probe set $\mathcal{P}=\{z_m\}_{m=1}^M\subset{\Omega}$. The \emph{absolute error} at a point $z\in{\Omega}$ is
\[
e(z)\;=\;\big|\,u(z)-u_h(z)\,\big|,
\]
and the \emph{root-mean-square error} (RMSE),
\[
\mathrm{RMSE}\;=\;\Bigg(\frac{1}{M}\sum_{m=1}^M \big|\,u(z_m)-u_h(z_m)\,\big|^2\Bigg)^{1/2}.
\]

\subsection*{Example 1: 1D Poisson, Dirichlet}
\begin{figure}[H]
\centering
 \begin{subfigure}{0.45\textwidth}
     \includegraphics[width=\textwidth]{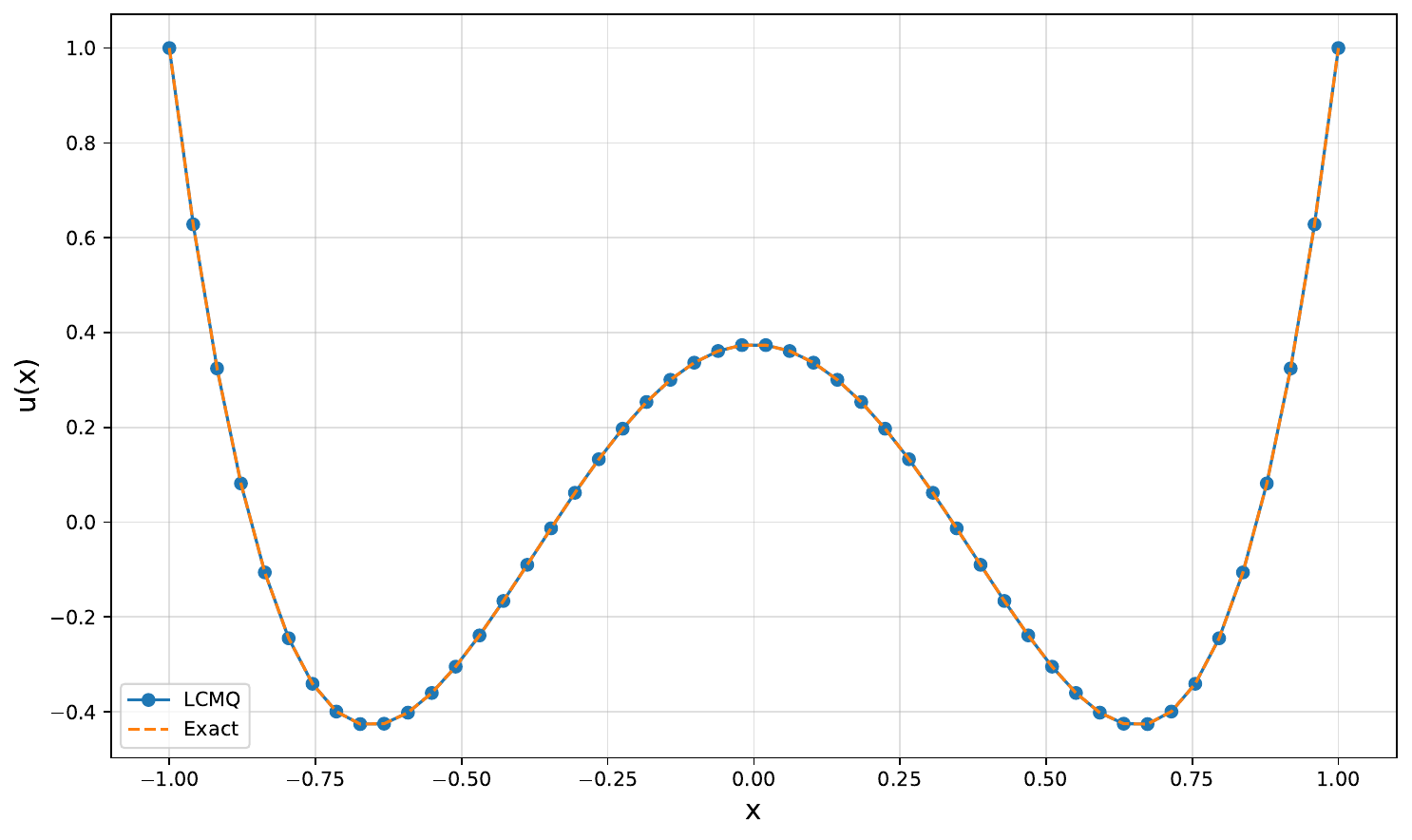}
     \caption{Numerical and exact solutions $u(x)$ on $[-1,1]$.}
     \label{fig:1a}
 \end{subfigure}
 \begin{subfigure}{0.45\textwidth}
     \includegraphics[width=\textwidth]{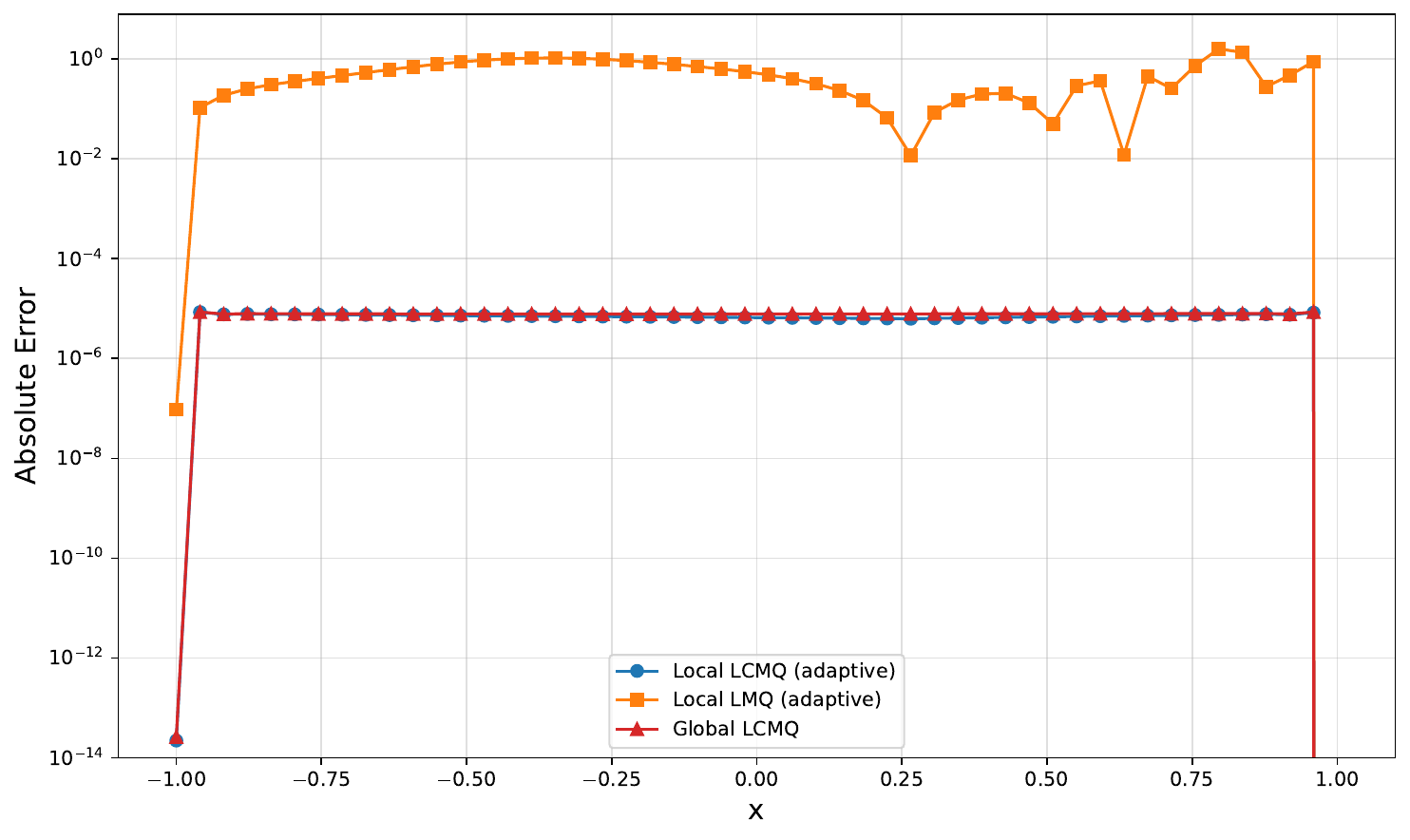}
     \caption{Absolute error $\lvert u_h-u^{\star}\rvert$ over $[-1,1]$.}
     \label{fig:1b}
 \end{subfigure}
  \begin{subfigure}{0.45\textwidth}
     \includegraphics[width=\textwidth]{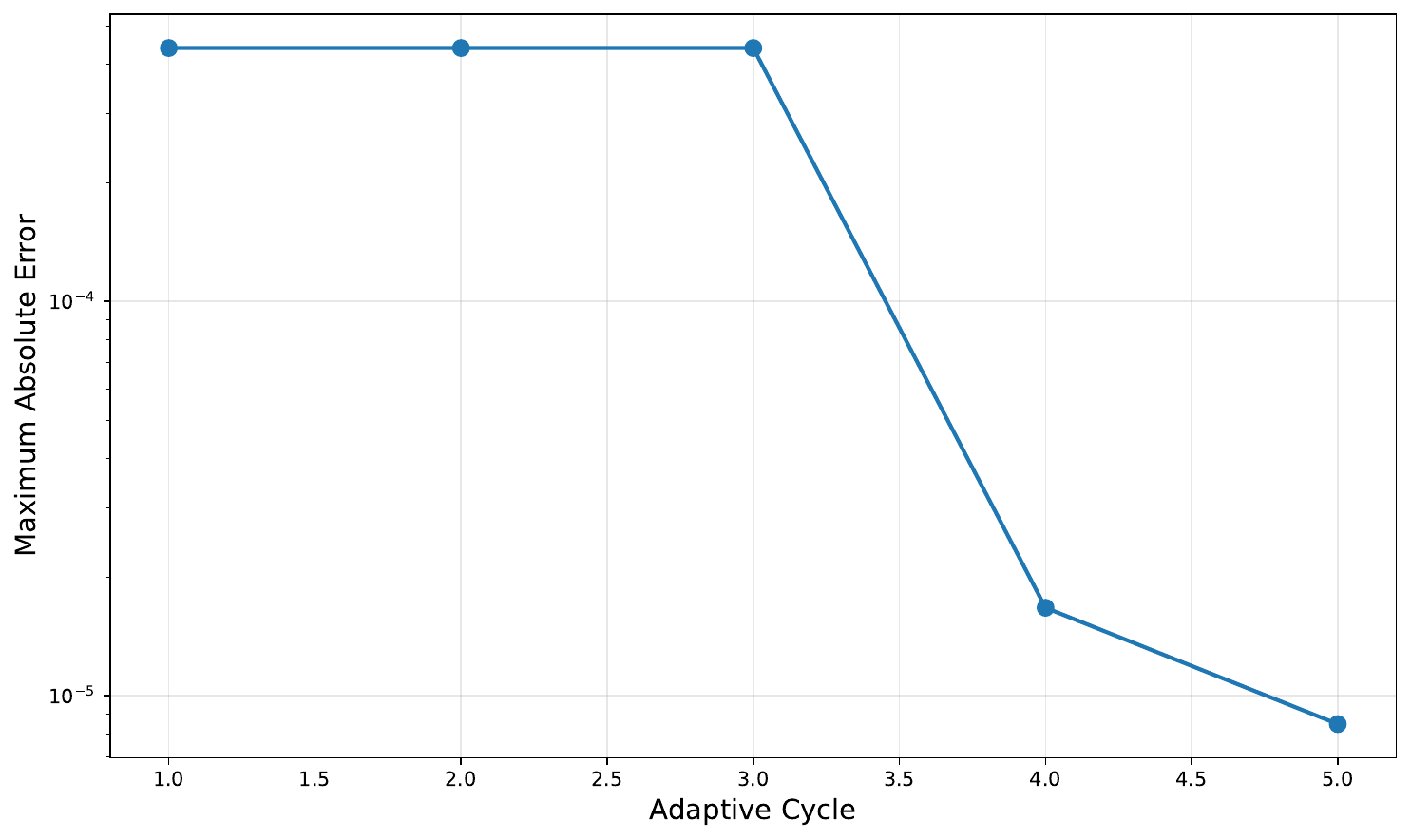}
     \caption{Maximum absolute error versus adaptive cycle.}
     \label{fig:1c}
 \end{subfigure}
 \begin{subfigure}{0.45\textwidth}
     \includegraphics[width=\textwidth]{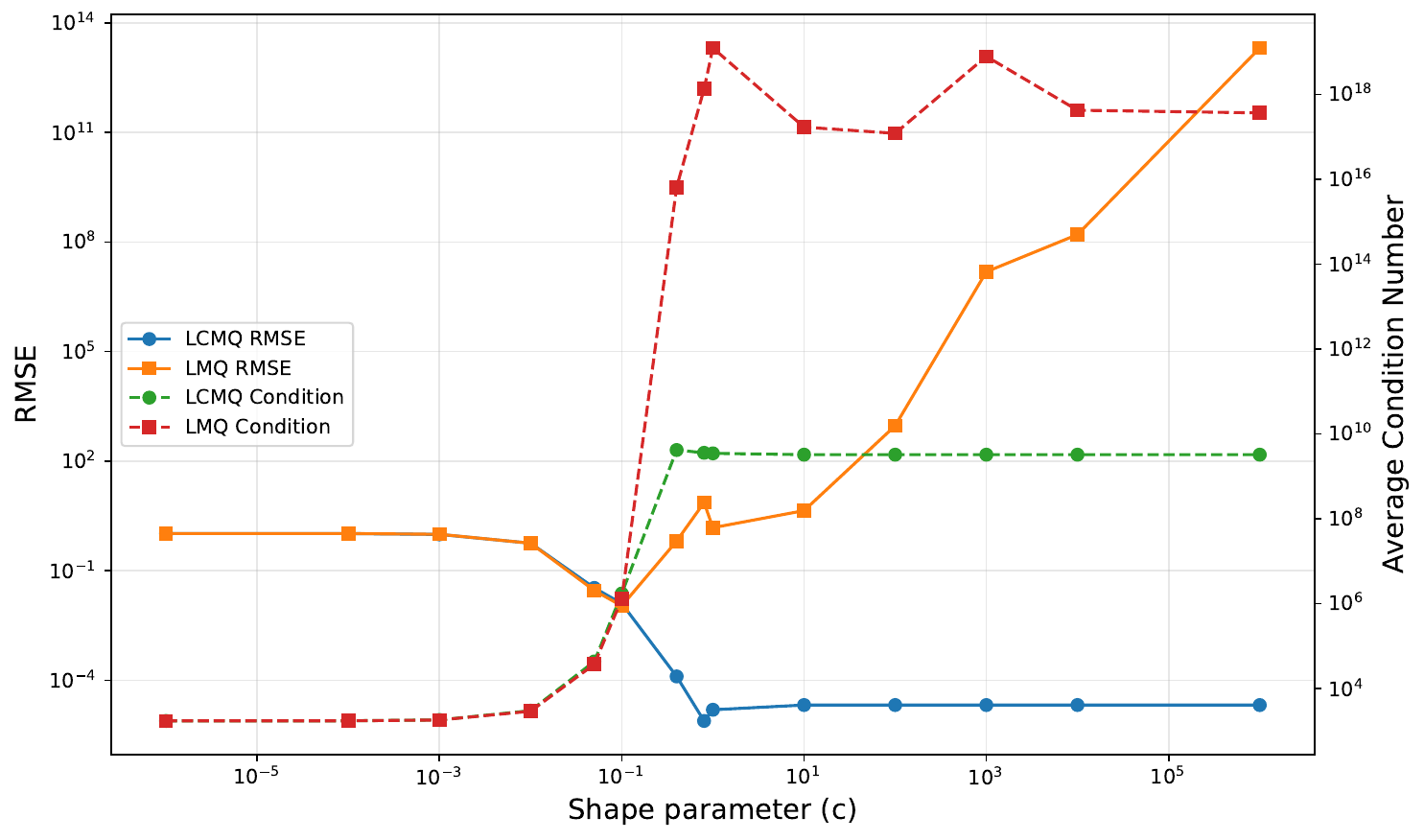}
     \caption{RMSE versus variation in the shape parameter $c$.}
     \label{fig:1d}
 \end{subfigure}
 \caption{Results for the 1D Poisson problem with adaptive $\mathcal{C}h$--LCMQ.}
 \label{fig:1d-results}
\end{figure}

\begin{figure}[H]
    \centering
    \includegraphics[width=0.5\linewidth]{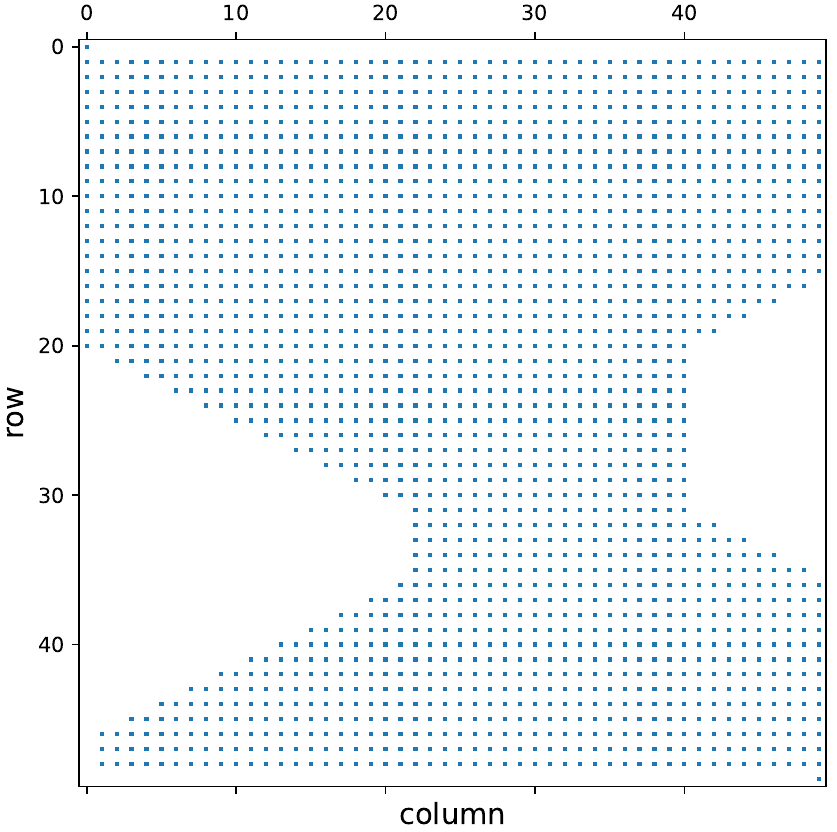}
    \caption{Cover size distribution for the 1D Poisson problem at the last adaptive cycle}
    \label{fig:1d_spy}
\end{figure}

Find $u:\,[-1,1]\to\mathbb{R}$ such that
\begin{equation}
u''(x)=\frac{105}{2}x^2-\frac{15}{2},\qquad u(-1)=u(1)=1.
\end{equation}
The exact solution is
\begin{equation}
u^\star(x)=\frac{35}{8}x^4-\frac{15}{4}x^2+\frac{3}{8}.
\end{equation}
In the 1D case, we apply the decision policy of Algorithm~\ref{alg:ACMQ} with the solver settings listed in Table~\ref{tab:solver-config-1d}. The adaptive procedure reduces the error over successive cycles (Figure~\ref{fig:1c}) and reaches the prescribed tolerance after five iterations by adjusting only the cover sizes, as summarized in Table~\ref{tab:1d_performance}. The final cover distribution (Figure~\ref{fig:1d_spy}) shows larger covers near the boundaries, where the error indicator $\mathcal{E}$ is highest~\cite{deng2023accuracy}. Compared with LMQ and global collocation, the adaptive $\mathcal{C}h$ method with LCMQ achieves higher accuracy (Figure~\ref{fig:1b}), while also exhibiting lower condition numbers and reduced sensitivity to the shape parameter (Figure~\ref{fig:1d}).

\begin{table}[!h]
\centering
\caption{Solver configuration for Examples 1 and 2}
\label{tab:solver-config-1d}
\begingroup
\renewcommand{\arraystretch}{1.15}
\begin{tabular}{@{}lll@{}}
\toprule
\textbf{Parameter (symbol)} & \textbf{Values in Example 1} & \textbf{Values in Example 2} \\ \midrule
Domain $[a,b]$                & $[-1,\,1]$                 & $[0,\,1], [0,\,1]$          \\
Nodes $N_{0}$                 & $50$                       & $20\times20$                \\
Cover size $\mathcal{C}_0$    & $7$                        & $50$                         \\ 
Shape parameter $c$           & $0.8$                      & $0.8$                       \\
Residual tol.\ $\tau$         & $1\times 10^{-5}$          & $8\times10^{-5}$            \\
\bottomrule
\end{tabular}
\endgroup
\end{table}

\begin{table}[!h]
\centering
\caption{Solver performance for the 1D Poisson problem on an Intel Core i7-6820HQ @ 2.70GHz}
\label{tab:1d_performance}
\begin{tabular}{@{}cccccccc@{}}
\toprule
Cycle & Nodes & Max Cover Size & RMSE & Absolute Max Error & Solve (s) & Error Estimation \& Adaptation (s)\\
\midrule
1 & 50 & 7 & 3.37e-04 & 4.39e-04 & 0.0071  & 0.0110\\
2 & 50 & 50 & 3.37e-04 & 4.39e-04 & 0.0084 & 0.0219\\
3 & 50 & 17 & 3.37e-04 & 4.39e-04 & 0.0070 & 0.0152\\
4 & 50 & 50 & 1.33e-05 & 1.67e-05 & 0.0131 & 0.0252\\
5 & 50 & 50 & 6.98e-06 & 8.49e-06 & 0.0121 & 0.0267\\
\bottomrule
\end{tabular}
\end{table}
\subsection*{Example 2: Poisson on the unit square, Dirichlet}
Now, let $\Omega=[0,1]^2$. it is required to find $u:\Omega\to\mathbb{R}$ such that
\begin{equation}
\nabla^2 u(x,y)=\sin(\pi x)\sin(\pi y),\qquad u|_{\partial\Omega}=0.
\end{equation}
The exact solution is
\begin{equation}
u^\star(x,y)=-\frac{1}{2\pi^2}\sin(\pi x)\sin(\pi y).
\end{equation}

\begin{figure}[!h]
\centering
     \includegraphics[width=0.7\textwidth]{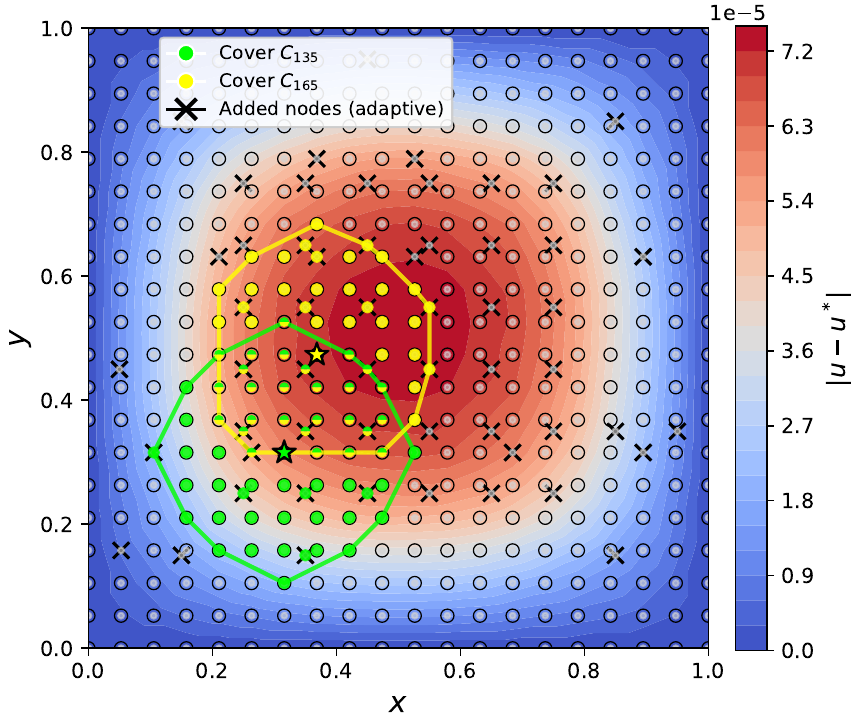}
     \caption{Numerical solution and absolute error of the 2D Poisson problem.}
\label{fig:2}
\end{figure}

\begin{figure}[!h]
\centering
\begin{subfigure}{0.7\textwidth}
    \includegraphics[width=1\textwidth]{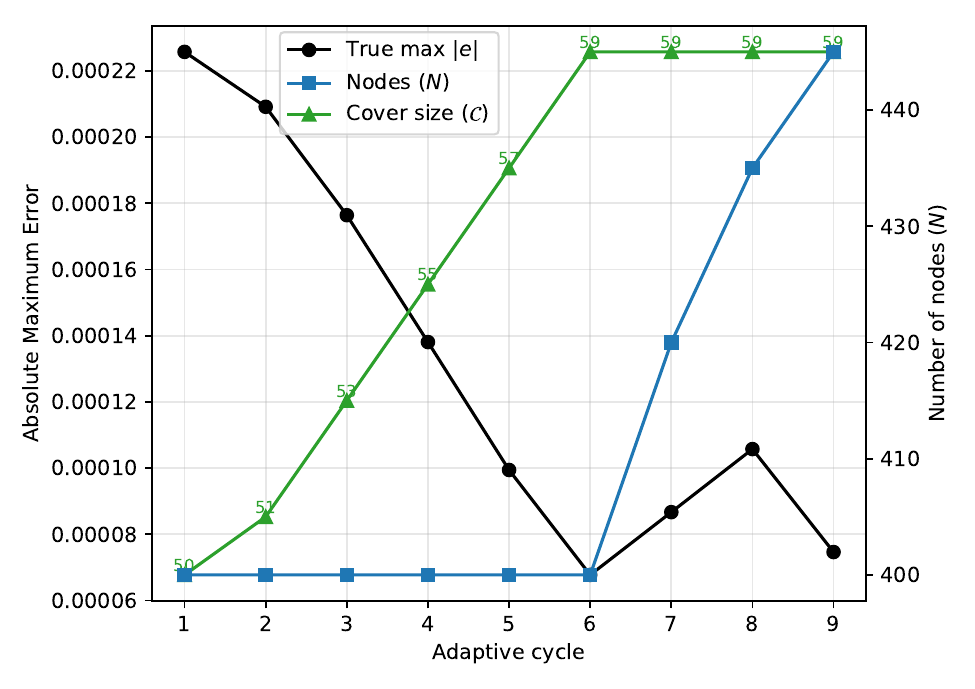}
    \caption{Absolute Max error vs Adaptive Cycle Number.}
    \label{fig:3a}
\end{subfigure}
\hfill
\begin{subfigure}{0.7\textwidth}
    \includegraphics[width=1\textwidth]{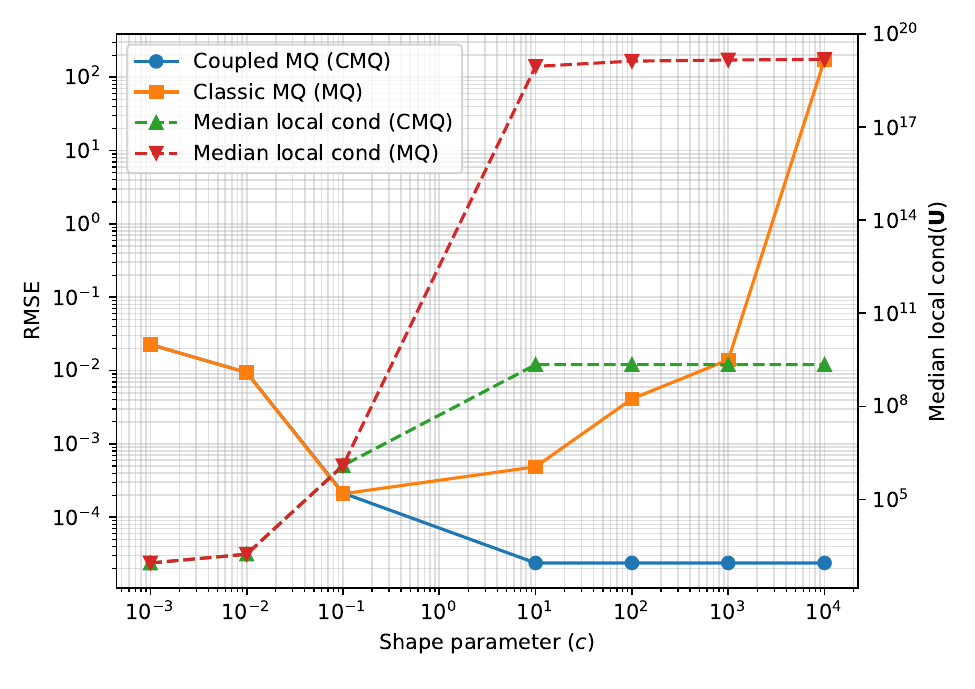}
    \caption{RMSE vs variation in the shape parameter.}
    \label{fig:3b}
\end{subfigure}
\caption{Error in the 2D Poisson problem.}
\label{Fig3}
\end{figure}

\begin{table}[h!]
\centering
\caption{Solver performance for the 2D Poisson problem on an Intel Core i7-6820HQ @ 2.70GHz}
\label{tab:2d_performance}
\begin{tabular}{@{}cccccccc@{}}
\toprule
Cycle & Nodes & Max Cover Size & RMSE & Absolute Max Error & Solve (s) & Error Estimation \& Adaptation (s)\\
\midrule
1  & 400 & 51 & 1.05e-04 & 2.26e-04  & 0.1013 & 0.0303 \\
2  & 400 & 53 & 9.70e-05 & 2.09e-04  & 0.0921 & 0.0295 \\
3  & 400 & 55 & 8.12e-05 & 1.76e-04  & 0.0944 & 0.0335 \\
4  & 400 & 57 & 6.29e-05 & 1.38e-04  & 0.1046 & 0.0342 \\
5  & 400 & 59 & 4.47e-05 & 9.94e-05  & 0.1023 & 0.0319 \\
6  & 420 & 59 & 3.08e-05 & 6.77e-05  & 0.1099 & 0.0342 \\
7  & 435 & 59 & 4.01e-05 & 8.67e-05  & 0.1217 & 0.0377 \\
8  & 445 & 59 & 5.22e-05 & 1.06e-04  & 0.1204 & 0.0373 \\
9  & 452 & 59 & 4.12e-05 & 7.46e-05  & 0.1196 & 0.0325 \\
\bottomrule
\end{tabular}
\end{table}

Adaptive $\mathcal{C}h$ with LCMQ achieves high accuracy (Figure~\ref{fig:2}), well below the tolerance in Table~\ref{tab:solver-config-1d}. The same figure illustrates two representative covers together with the nodes inserted by the adaptive algorithm. The solver terminates after 9 cycles with RMSE on the order of $10^{-5}$ (Figure~\ref{fig:3a}), performing cover-size refinement first and local node insertion thereafter, as summarized in Table~\ref{tab:2d_performance}.  For brevity, the table reports the maximum cover size and total number of nodes for selected cycles. It is noteworthy that in 1D the error–estimation/adaptation stage exceeded the solve time, whereas in 2D it fell below it, because the error estimation work grows roughly linearly with cover size while the banded assembly/linear solve escalates superlinearly with the larger 2D problem size and bandwidth. Finally,  LCMQ demonstrates greater stability than MQ  over a longer range or shape parameter values (Figure~\ref{fig:3b}).

\section{Conclusions}
We have extended the global CMQ method to a local one within an adaptive framework that relies exclusively on node clouds and local covers, yielding a truly meshless approximation scheme. The Adaptive $\mathcal{C}h$ Method with LCMQ inherits the desirable shape-parameter low sensitivity of CMQs while operating with sparse, banded matrices assembled from local covers. The adaptive logic, which prioritizes local $\mathcal{C}$-enrichment and resorts to $h$-refinement only when necessary, is simple to implement yet capable of targeting regions where the PDE operator is under-resolved.

The numerical experiments for 1D and 2D Poisson problems demonstrate that the method achieves accuracy comparable to global collocation while avoiding the dense matrices and severe conditioning issues typically associated with global RBF approaches. In the 1D example, a fixed set of $50$ nodes suffices to reach a maximum absolute error below $10^{-5}$ by adjusting only the cover sizes. In the 2D unit-square problem, the algorithm refines both covers and nodes, reducing the RMSE to the order of $10^{-5}$ with fewer than $500$ nodes and cover sizes capped at $\mathcal{C}_{\max}=59$. In both cases the residual-based a~posteriori indicator successfully identifies regions requiring refinement, and the classifier $s_g$ effectively biases the method toward $\mathcal{C}$-enrichment or $h$-refinement.

From a computational standpoint, the experiments indicate that the additional work introduced by the adaptive engine remains moderate compared with the cost of the linear solves. In 1D, error estimation and adaptation cost slightly more than the solve, but the absolute times are small because the global system is only $50\times 50$. In 2D, the situation reverses: the linear solve dominates the per-cycle cost, while the indicator evaluation scales essentially linearly with the number and size of covers. The global matrix assembled from local LCMQ kernels stays banded, and its bandwidth is directly related to the maximum cover size, which is controlled by the algorithmic parameters. This structure is advantageous for direct solvers and also for preconditioned iterative methods.

The reported shape-parameter sweeps confirm that LCMQ-based collocation is significantly less sensitive to the choice of shape parameter than standard MQ-based local methods, both in 1D and in 2D. This removes the need for problem-specific tuning of $c$ in the examples considered and suggests that a single moderate value of $c$ can be used across a family of problems without substantial loss of accuracy. This property is particularly useful in adaptive settings, where the local fill distance and cover geometry change throughout the computation.

Benchmarks suggest that LCMQ collocation within the Adaptive $\mathcal{C}h$ framework offers a reliable, scalable pathway to solving PDEs with reduced sensitivity to $c$. Extensions of this work include: (i) investigating more sophisticated marking strategies and convergence guarantees for the residual-based indicator; (ii) studying preconditioning and parallelization strategies that exploit the banded structure and the locality of covers for large-scale three-dimensional problems. These directions build directly on the components introduced here and can be addressed within the same covers-and-nodes framework.

\subsection*{Declaration of generative AI and AI-assisted technologies in the writing process}

During the preparation of this work the author(s) used "Overleaf embedded AI Assistance and Writeful" in order to proofread and ensure language coherence across the article. After using this tool/service, the author(s) reviewed and edited the content as needed and take(s) full responsibility for the content of the publication.

\subsection*{Acknowledgements}
This work was self-funded by the author. No funding agencies or grants supported this work. 

\subsection*{Author contributions}
Material preparation, data collection and analysis were performed by Ahmed E. Seleit.

\subsection*{Financial disclosure}
None.

\subsection*{Conflict of interest}
The author declares no potential conflict of interests.

\subsection*{Data Availability}
No data was used for the research described in this article.

\bibliographystyle{ieeetr}
\bibliography{references}
\end{document}